\documentclass{article}
\usepackage{amssymb}
\usepackage{graphicx}
\usepackage{amsmath}

\input{tcilatex}
\begin{document}

\begin{center}
\bigskip {\LARGE A note on three dimensional good sets}

\bigskip \bigskip{\large K. GOWRI NAVADA}
\end{center}

\bigskip

Department of Mathematics, Periyar University, Salem - 636011, India

E-mail: gnavada@yahoo.com

\textbf{2000 Mathematics Subject Classification:} \quad primary 60A05,
47A35, \ secondary 28D05, 37Axx

\textbf{Abstract:.} We show that as in the case of n- fold
Cartesian product for $n\geq 4$, even in $3$-fold Cartesian
product, a related component need not be full component.

\textbf{Key words.}\quad Good set; full set; full component; related
component; geodesic; boundary of a good set. \

\textbf{Introduction and Preliminaries} The purpose of this note
is to answer two questions about good sets
raised in [3] and [4] for the case $n=3.$ \\

Let $X_{1},X_{2},\ldots ,X_{n}$ be nonempty sets and let $\Omega
=X_{1}\times X_{2}\times \cdots \times X_{n}$ be their Cartesian product. We
will write $\overset{\rightarrow}{x}$ to denote a point $(x_{1},x_{2},\ldots
,x_{n})\in \Omega .$

For each $1\leq i\leq n,\ \Pi_{i}$ denotes the canonical projection of $%
\Omega $ onto $X_{i}.$\newline

A subset $S\subset \Omega $ is said to be \textit{good}, if every complex
valued function $f$ on $S$ is of the form:
\begin{equation}
f(x_{1},x_{2},\ldots ,x_{n})=u_{1}(x_{1})+u_{2}(x_{2})+\cdots
+u_{n}(x_{n}),\ (x_{1},x_{2},\ldots ,x_{n})\in S,  \tag{1}
\end{equation}%
\qquad \qquad for suitable functions $u_{1},u_{2},\ldots ,u_{n}$ on $%
X_{1},X_{2},\ldots ,X_{n}$ respectively ([3], p. 181).\newline

For a good set $S$, a subset $B\subset \bigcup_{i=1}^{n}\Pi _{i}S$ is said
to be a \textit{boundary set of} $S$, if for any complex valued function $U$
on $B$ and for any $f:S\longrightarrow \mathbb{C}$ the equation $(1)$
subject to
\begin{equation*}
u_{i}|_{B\cap \Pi _{i}S}=U|_{B\cap \Pi _{i}S},\ 1\leq i\leq n,
\end{equation*}%
admits a unique solution. For a good set there always exists a boundary set
([3], p. 187). \newline

A subset $S\subset \Omega $ is said to be \textit{full}, if $S$ is maximal
good set in $\Pi_{1}S\times \Pi_{2}S\times\cdots \times\Pi _{n}S$.

A set $S \subset \Omega$ is full if and only if it has a boundary consisting
of $n-1$ points ([3], Theorem 3, page 185).\newline

If a set $S$ is good, maximal full subsets of $S$ form a partition of $S$.
They are called \textit{full components} of $S$ ([3], p. 183).\newline

Two points $\overset{\rightarrow}{x},\overset{\rightarrow}{y}$ in a good set
$S$ are said to be \textit{related}, denoted by $\overset{\rightarrow}{x} R
\overset{\rightarrow}{y}$, if there exists a finite subset of $S$ which is
full and contains both $\overset{\rightarrow}{x}$ and $\overset{\rightarrow}{%
y}$. $R$ is an equivalence relation, whose equivalence classes are called
\textit{related components} of $S$. The related components of $S$ are full
subsets of $S$ (ref. [3]).

\bigskip First we prove that when the dimension $n=3$, a full component need
not be a related component, by giving an example of a full set with
infinitely many related components.

Consider a countable set $T$ which consists of the following points:

$\overset{\rightarrow}{a_1}=(x_{1},x_{2},x_{3})$

$\overset{\rightarrow}{a_2}=(y_{1},y_{2},x_{3})$

$\overset{\rightarrow}{a_3}=(y_{1,}x_{2},z_{3})$ \medskip

$\overset{\rightarrow}{a_4}=(\alpha _{1},\alpha _{2},\alpha _{3})$

$\overset{\rightarrow}{a_5}=(\alpha _{4},\alpha _{5},\alpha _{3})$

$\overset{\rightarrow}{a_6}=(\alpha _{1},\alpha _{5},z_{3})$

$\overset{\rightarrow}{a_7}=(\alpha _{4},\alpha _{2},x_{3})$

$\overset{\rightarrow}{a_8}=(x_{1},y_{2},\alpha _{3})$ \medskip

$\overset{\rightarrow}{a_9}=(\alpha _{6},\alpha _{7},\alpha _{8})$

$\overset{\rightarrow}{a_{10}}=(\alpha _{9},\alpha _{10},\alpha _{8})$

$\overset{\rightarrow}{a_{11}}=(\alpha _{6},\alpha _{10},\alpha_{3})$

$\overset{\rightarrow}{a_{12}}=(\alpha _{9},\alpha _{7},x_{3})$

$\overset{\rightarrow}{a_{13}}=(x_{1},\alpha _{2},\alpha _{8})$ \medskip

$.....$

$......$

$\overset{\rightarrow}{a}_{5n-1}=(\alpha _{5n-4},\alpha _{5n-3},\alpha
_{5n-2})$

$\overset{\rightarrow}{a}_{5n}=(\alpha _{5n-1},\alpha _{5n},\alpha _{5n-2})$

$\overset{\rightarrow}{a}_{5n+1}=(\alpha _{5n-4},\alpha _{5n},\alpha
_{5n-7}) $

$\overset{\rightarrow}{a}_{5n+2}=(\alpha _{5n-1},\alpha _{5n-3},x_{3})$

$\overset{\rightarrow}{a}_{5n+3}=(x_{1},\alpha _{5n-8},\alpha _{5n-2})$

$....$

$......$

Call the first three points of $T$ as $D_{0}$ and for $n\geqslant 1,$ let $%
D_{n}$ denote the first $3+5n$ points of $S.$ Let $A_{0}=D_{0}$ and for $%
n\geq 1$ let $A_{n}=D_{n}\setminus D_{n-1}$. Then it is easy to see that
every $D_{n}$ is good and has three point boundary. All the three points of
the boundary of $D_{n}$ cannot come from the coordinates of points in $%
D_{n-1}$: because, if all of them occur as coordinates in $D_{n-1},$ they
form a boundary for $D_{n-1}$. Given any function $f$ on $D_{n},$ there is a
solution $u_{1},u_{2},u_{3}$ on $D_{n-1}$ such that
\begin{equation*}
f(w_{1},w_{2},w_{3})=u_{1}(w_{1})+u_{2}(w_{2})+u_{3}(w_{3}),\
(w_{1},w_{2},w_{3})\in D_{n-1}.
\end{equation*}%
But then $f(\overset{\rightarrow }{a}_{5n+3})$ fixes the value of $%
u_{3}(\alpha _{5n-2})$ by the following equation:
\begin{equation*}
u_{3}(\alpha _{5n-2})=f(\overset{\rightarrow }{a}%
_{5n+3})-u_{1}(x_{1})-u_{2}(\alpha _{5n-8})
\end{equation*}%
When we substitute this value of $u_{3}(\alpha _{5n-2})$ in the remaining
four points of $A_{n}$, we get a set of linearly dependent equations. This
shows that the boundary of $D_{n}$ contains at least one of the five
coordinates, $\alpha _{5n-4},\alpha _{5n-3},\alpha _{5n-2},\alpha _{5n-1}$
or $\alpha _{5n},$ which are introduced in $A_{n}.$ One can observe the
following properties of the points in the set $A_{n}$: any $k$ points of $%
A_{n}$ has at least $k$ coordinates introduced in $A_{n}$. (i.e,
they do not
occur as coordinates in $D_{n-1}$)$.$ If we take a singleton $\{\overset{%
\rightarrow }{a_{i}}\}$ in $D_{n-1}$, any set of $k$ points of $A_{n}$ has
at least $(k+1)$ coordinates which do not occur as coordinates of $\overset{%
\rightarrow }{a_{i}}.$

$T$ is good as every finite subset of $T$ is good. It cannot have a boundary
$B$ with more than two points: If $|B|=3,$ we can choose a $n$ sufficiently
large such that all the three points of $b$ occur as coordinates in $%
D_{n-1}. $ Then $B$ is a boundary of $D_{n}$ which is not possible as
observed above. If $|B|>3,$ we can choose $n$ sufficiently large so that $%
k=|B\cap \cup _{i=1}^{3}\Pi _{i}D_{n}\geqslant 4.$ Then these $k$ points
form a boundary of $D_{n}$ which is again not possible. So the boundary of $%
T $ consists of only two points which shows that $T$ is full.

We prove that no finite subset $A$ of $T$ other than singleton is full: Set $%
|A\cap A_{i}|=k_{i}$ for $i\geq 0$. Let $i_{1}<i_{2}<\cdots <i_{l}$ be such
that $k_{i_{j}}\neq 0$ for $j=1,2,\ldots l$ and $k_{i}=0$ for all other $i$.
If $k_{i_{1}}>1$, as no subset other than singleton of $A_{n}$ is full, the
set $A\cap A_{i_{1}}$ is not full. When we add the points of $A\cap
A_{i_{2}} $ to $A\cap A_{i_{1}}$ ( as we are adding $k_{i_{2}}$ points) we
will be adding at least $k_{i_{2}}$ new coordinates. So the set $A\cap
(A_{i_{1}}\cup A_{i_{2}})$ is not full. Similarly when we keep adding $A\cap
A_{i_{j}}$ to the set $A\cap (\cup _{k<j}A_{i_{k}})$ the number of
coordinates added is at least equal to the number of points added. So at
each step $A\cap (\cup _{k\leq j}A_{i_{k}})$ is not full. In this way we get
$A=A\cap (\cup _{k\leq l}A_{i_{k}})$ is also not full. If $k_{i_{1}}=1$, in
the first step when we add points of $A\cap A_{i_{2}}$ to the singleton set $%
A\cap A_{i_{1}}$ the new coordinates added is at least $k_{i_{2}}+1$. So $%
A\cap (A_{i_{1}}\cup A_{i_{2}})$ is not full. In the remaining steps as we
keep adding points from $A\cap A_{i_{j}}$, the number of coordinates added
is at least equal to the number of points added. So in the end we get $A$ is
not full.

For any $n,$ let $\overset{\rightarrow }{b_{n}}=(\alpha _{5n-1},y_{2},z_{3})$
and consider the set $F_{n}=D_{n}\cup \overset{\rightarrow }{b_{n}}$. We
show that the geodesic between the points $\overset{\rightarrow }{a_{1}}$
and $\overset{\rightarrow }{a}_{5n+3}$ in $F_{n}$ is the whole set $F_{n}$.
To show that $F_{n}$ is full, consider the matrix $M_{n}$ whose rows
correspond to the points $\overset{\rightarrow }{a_{2}},\overset{\rightarrow
}{a_{3}},\cdots ,\overset{\rightarrow }{a}_{5n+3},\overset{\rightarrow }{%
b_{n}}$ and columns correspond to the coordinates $y_{1},y_{2},z_{3},\alpha
_{1},\alpha _{2},...,\alpha _{5n}.$This is a $5n+3\times 5n+3$ matrix:

$%
\begin{bmatrix}
1 & 1 & 0 & 0 & 0 & 0 & 0 & 0 & . & . & 0 & 0 & 0 & 0 & 0 \\
1 & 0 & 1 & 0 & 0 & 0 & 0 & 0 & . & . & 0 & 0 & 0 & 0 & 0 \\
0 & 0 & 0 & 1 & 1 & 1 & 0 & 0 & . & . & 0 & 0 & 0 & 0 & 0 \\
0 & 0 & 0 & 0 & 0 & 1 & 1 & 1 & . & . & 0 & 0 & 0 & 0 & 0 \\
0 & 0 & 1 & 1 & 0 & 0 & 0 & 1 & . & . & 0 & 0 & 0 & 0 & 0 \\
0 & 0 & 0 & 0 & 1 & 0 & 1 & 0 & . & . & 0 & 0 & 0 & 0 & 0 \\
0 & 1 & 0 & 0 & 0 & 1 & 0 & 0 & . & . & 0 & 0 & 0 & 0 & 0 \\
0 & 0 & 0 & 0 & 0 & 0 & 0 & 0 & . & . & 1 & 1 & 1 & 0 & 0 \\
. & . & . & . & . & . & . & . & . & . & . & . & . & . & . \\
. & . & . & . & . & . & . & . & . & . & . & . & . & . & . \\
0 & 0 & 0 & 0 & 0 & 0 & 0 & 0 & . & . & 0 & 0 & 1 & 1 & 1 \\
0 & 0 & 0 & 0 & 0 & 1 & 0 & 0 & . & . & 1 & 0 & 0 & 0 & 1 \\
0 & 0 & 0 & 0 & 0 & 0 & 0 & 0 & . & . & 0 & 1 & 0 & 1 & 0 \\
0 & 0 & 0 & 0 & 1 & 0 & 0 & 0 & . & . & 0 & 0 & 1 & 0 & 0 \\
0 & 1 & 1 & 0 & 0 & 0 & 0 & 0 & . & . & 0 & 0 & 0 & 1 & 0%
\end{bmatrix}%
$

It has an inverse given by

$
{\tiny
\begin{bmatrix}
\frac{2}{3}\hspace{-0.15in} & \frac{1}{3}\hspace{-0.15in} & \frac{1}{3}\hspace{%
-0.2in} & \frac{1}{3}\hspace{-0.15in} & \frac{-1}{3}\hspace{-0.15in} & \frac{-1%
}{3}\hspace{-0.15in} & \frac{-2}{3}\hspace{-0.15in} & 0\hspace{-0.15in} & 0%
\hspace{-0.15in} & 0\hspace{-0.15in} & .\hspace{-0.15in} & .\hspace{-0.15in} & .%
\hspace{-0.15in} & .\hspace{-0.15in} & .\hspace{-0.15in} & .\hspace{-0.15in} & .%
\hspace{-0.15in} & .\hspace{-0.15in} & .\hspace{-0.15in} & .\hspace{-0.15in} & .%
\hspace{-0.15in} & .\hspace{-0.15in} & .\hspace{-0.15in} & .\hspace{-0.15in} & .%
\hspace{-0.15in} & .\hspace{-0.15in} & .\hspace{-0.15in} & .\hspace{-0.15in} & .%
\hspace{-0.15in} & 0 \\
\hspace{-0.15in} & \hspace{-0.15in} & \hspace{-0.15in} &
\hspace{-0.15in} & \hspace{-0.15in} & \hspace{-0.15in} &
\hspace{-0.15in} & \hspace{-0.15in} & \hspace{-0.15in} &
\hspace{-0.15in} & \hspace{-0.15in} & \hspace{-0.15in} &
\hspace{-0.15in} & \hspace{-0.15in} & \hspace{-0.15in} &
\hspace{-0.15in} & \hspace{-0.15in} & \hspace{-0.15in} &
\hspace{-0.15in} & \hspace{-0.15in} & \hspace{-0.15in} &
\hspace{-0.15in} & \hspace{-0.15in} & \hspace{-0.15in} &
\hspace{-0.15in} & \hspace{-0.15in} & \hspace{-0.15in} &
\hspace{-0.15in} &
\hspace{-0.15in} &  \\
\frac{1}{3}\hspace{-0.15in} & \frac{-1}{3}\hspace{-0.15in} & \frac{-1}{3}%
\hspace{-0.15in} & \frac{-1}{3}\hspace{-0.15in} &
\frac{1}{3}\hspace{-0.15in} & \frac{1}{3}\hspace{-0.15in} &
\frac{2}{3}\hspace{-0.15in} & 0\hspace{-0.15in} &
0\hspace{-0.15in} & 0\hspace{-0.15in} & .\hspace{-0.15in} & .\hspace{-0.15in} & .%
\hspace{-0.15in} & .\hspace{-0.15in} & .\hspace{-0.15in} & .\hspace{-0.15in} & .%
\hspace{-0.15in} & .\hspace{-0.15in} & .\hspace{-0.15in} & .\hspace{-0.15in} & .%
\hspace{-0.15in} & .\hspace{-0.15in} & .\hspace{-0.15in} & .\hspace{-0.15in} & .%
\hspace{-0.15in} & .\hspace{-0.15in} & ..\hspace{-0.15in} & .\hspace{-0.15in} & .%
\hspace{-0.15in} & 0 \\
\hspace{-0.15in} & \hspace{-0.15in} & \hspace{-0.15in} &
\hspace{-0.15in} & \hspace{-0.15in} & \hspace{-0.15in} &
\hspace{-0.15in} & \hspace{-0.15in} & \hspace{-0.15in} &
\hspace{-0.15in} & \hspace{-0.15in} & \hspace{-0.15in} &
\hspace{-0.15in} & \hspace{-0.15in} & \hspace{-0.15in} &
\hspace{-0.15in} & \hspace{-0.15in} & \hspace{-0.15in} &
\hspace{-0.15in} & \hspace{-0.15in} & \hspace{-0.15in} &
\hspace{-0.15in} & \hspace{-0.15in} & \hspace{-0.15in} &
\hspace{-0.15in} & \hspace{-0.15in} & \hspace{-0.15in} &
\hspace{-0.15in} &
\hspace{-0.15in} &  \\
\frac{-2}{3}\hspace{-0.15in} & \frac{2}{3}\hspace{-0.15in} & \frac{-1}{3}%
\hspace{-0.15in} & \frac{-1}{3}\hspace{-0.15in} &
\frac{1}{3}\hspace{-0.15in} & \frac{1}{3}\hspace{-0.15in} &
\frac{2}{3}\hspace{-0.15in} & 0\hspace{-0.15in} &
0\hspace{-0.15in} & 0\hspace{-0.15in} & .\hspace{-0.15in} & .\hspace{-0.15in} & .%
\hspace{-0.15in} & .\hspace{-0.15in} & .\hspace{-0.15in} & .\hspace{-0.15in} & .%
\hspace{-0.15in} & .\hspace{-0.15in} & .\hspace{-0.15in} & .\hspace{-0.15in} & .%
\hspace{-0.15in} & .\hspace{-0.15in} & .\hspace{-0.15in} & .\hspace{-0.15in} & .%
\hspace{-0.15in} & .\hspace{-0.15in} & .\hspace{-0.15in} & .\hspace{-0.15in} & .%
\hspace{-0.15in} & 0 \\
\hspace{-0.15in} & \hspace{-0.15in} & \hspace{-0.15in} &
\hspace{-0.15in} & \hspace{-0.15in} & \hspace{-0.15in} &
\hspace{-0.15in} & \hspace{-0.15in} & \hspace{-0.15in} &
\hspace{-0.15in} & \hspace{-0.15in} & \hspace{-0.15in} &
\hspace{-0.15in} & \hspace{-0.15in} & \hspace{-0.15in} &
\hspace{-0.15in} & \hspace{-0.15in} & \hspace{-0.15in} &
\hspace{-0.15in} & \hspace{-0.15in} & \hspace{-0.15in} &
\hspace{-0.15in} & \hspace{-0.15in} & \hspace{-0.15in} &
\hspace{-0.15in} & \hspace{-0.15in} & \hspace{-0.15in} &
\hspace{-0.15in} &
\hspace{-0.15in} &  \\
\frac{1}{6}\hspace{-0.15in} & \frac{1}{6}\hspace{-0.15in} & \frac{5}{6}\hspace{%
-0.2in} & \frac{-1}{6}\hspace{-0.15in} & \frac{1}{6}\hspace{-0.15in} & \frac{1%
}{6}\hspace{-0.15in} & \frac{-1}{6}\hspace{-0.15in} & \frac{1}{2}\hspace{-0.2in%
} & \frac{1}{2}\hspace{-0.15in} & \frac{-1}{2}\hspace{-0.15in} & \frac{-1}{2}%
\hspace{-0.15in} & -1\hspace{-0.15in} & 0\hspace{-0.15in} & 0\hspace{-0.15in} & 0%
\hspace{-0.15in} & .\hspace{-0.15in} & .\hspace{-0.15in} & .\hspace{-0.15in} & .%
\hspace{-0.15in} & .\hspace{-0.15in} & .\hspace{-0.15in} & .\hspace{-0.15in} & .%
\hspace{-0.15in} & .\hspace{-0.15in} & .\hspace{-0.15in} & .\hspace{-0.15in} & .%
\hspace{-0.15in} & .\hspace{-0.15in} & .\hspace{-0.15in} & 0 \\
\hspace{-0.15in} & \hspace{-0.15in} & \hspace{-0.15in} &
\hspace{-0.15in} & \hspace{-0.15in} & \hspace{-0.15in} &
\hspace{-0.15in} & \hspace{-0.15in} & \hspace{-0.15in} &
\hspace{-0.15in} & \hspace{-0.15in} & \hspace{-0.15in} &
\hspace{-0.15in} & \hspace{-0.15in} & \hspace{-0.15in} &
\hspace{-0.15in} & \hspace{-0.15in} & \hspace{-0.15in} &
\hspace{-0.15in} & \hspace{-0.15in} & \hspace{-0.15in} &
\hspace{-0.15in} & \hspace{-0.15in} & \hspace{-0.15in} &
\hspace{-0.15in} & \hspace{-0.15in} & \hspace{-0.15in} &
\hspace{-0.15in} &
\hspace{-0.15in} &  \\
\frac{1}{6}\hspace{-0.15in} & \frac{-1}{6}\hspace{-0.15in} & \frac{-1}{6}%
\hspace{-0.15in} & \frac{-1}{6}\hspace{-0.15in} &
\frac{1}{6}\hspace{-0.15in} &
\frac{1}{6}\hspace{-0.15in} & \frac{-1}{6}\hspace{-0.15in} & \frac{-1}{2}%
\hspace{-0.15in} & \frac{-1}{2}\hspace{-0.15in} &
\frac{1}{2}\hspace{-0.15in} &
\frac{1}{2}\hspace{-0.15in} & 1\hspace{-0.15in} & 0\hspace{-0.15in} & 0\hspace{%
-0.2in} & 0\hspace{-0.15in} & .\hspace{-0.15in} & .\hspace{-0.15in} & .\hspace{%
-0.2in} & .\hspace{-0.15in} & .\hspace{-0.15in} & .\hspace{-0.15in} & .\hspace{%
-0.2in} & .\hspace{-0.15in} & .\hspace{-0.15in} & .\hspace{-0.15in} & .\hspace{%
-0.2in} & .\hspace{-0.15in} & .\hspace{-0.15in} & .\hspace{-0.15in} & 0 \\
\hspace{-0.15in} & \hspace{-0.15in} & \hspace{-0.15in} &
\hspace{-0.15in} & \hspace{-0.15in} & \hspace{-0.15in} &
\hspace{-0.15in} & \hspace{-0.15in} & \hspace{-0.15in} &
\hspace{-0.15in} & \hspace{-0.15in} & \hspace{-0.15in} &
\hspace{-0.15in} & \hspace{-0.15in} & \hspace{-0.15in} &
\hspace{-0.15in} & \hspace{-0.15in} & \hspace{-0.15in} &
\hspace{-0.15in} & \hspace{-0.15in} & \hspace{-0.15in} &
\hspace{-0.15in} & \hspace{-0.15in} & \hspace{-0.15in} &
\hspace{-0.15in} & \hspace{-0.15in} & \hspace{-0.15in} &
\hspace{-0.15in} &
\hspace{-0.15in} &  \\
\frac{-1}{3}\hspace{-0.15in} & \frac{1}{3}\hspace{-0.15in} & \frac{1}{3}%
\hspace{-0.15in} & \frac{1}{3}\hspace{-0.15in} &
\frac{-1}{3}\hspace{-0.15in} & \frac{-1}{3}\hspace{-0.15in} &
\frac{1}{3}\hspace{-0.15in} & 0\hspace{-0.15in} &
0\hspace{-0.15in} & 0\hspace{-0.15in} & 0\hspace{-0.15in} &
0\hspace{-0.15in} & 0\hspace{-0.15in} & 0\hspace{-0.15in} &
0\hspace{-0.15in} & .\hspace{-0.15in} & .\hspace{-0.15in} &
.\hspace{-0.15in} & .\hspace{-0.15in} & .\hspace{-0.15in} &
.\hspace{-0.15in} & .\hspace{-0.15in} & .\hspace{-0.15in} &
.\hspace{-0.15in} & .\hspace{-0.15in} & .\hspace{-0.15in} &
.\hspace{-0.15in} & .\hspace{-0.15in}
& .\hspace{-0.15in} & 0 \\
\hspace{-0.15in} & \hspace{-0.15in} & \hspace{-0.15in} &
\hspace{-0.15in} & \hspace{-0.15in} & \hspace{-0.15in} &
\hspace{-0.15in} & \hspace{-0.15in} & \hspace{-0.15in} &
\hspace{-0.15in} & \hspace{-0.15in} & \hspace{-0.15in} &
\hspace{-0.15in} & \hspace{-0.15in} & \hspace{-0.15in} &
\hspace{-0.15in} & \hspace{-0.15in} & \hspace{-0.15in} &
\hspace{-0.15in} & \hspace{-0.15in} & \hspace{-0.15in} &
\hspace{-0.15in} & \hspace{-0.15in} & \hspace{-0.15in} &
\hspace{-0.15in} & \hspace{-0.15in} & \hspace{-0.15in} &
\hspace{-0.15in} &
\hspace{-0.15in} &  \\
\frac{-1}{6}\hspace{-0.15in} & \frac{1}{6}\hspace{-0.15in} & \frac{1}{6}%
\hspace{-0.15in} & \frac{1}{6}\hspace{-0.15in} &
\frac{-1}{6}\hspace{-0.15in} &
\frac{5}{6}\hspace{-0.15in} & \frac{1}{6}\hspace{-0.15in} & \frac{1}{2}\hspace{%
-0.2in} & \frac{1}{2}\hspace{-0.15in} & \frac{-1}{2}\hspace{-0.15in} & \frac{-1%
}{2}\hspace{-0.15in} & -1\hspace{-0.15in} & 0\hspace{-0.15in} &
0\hspace{-0.15in} & 0\hspace{-0.15in} & .\hspace{-0.15in} &
.\hspace{-0.15in} & .\hspace{-0.15in} & .\hspace{-0.15in} &
.\hspace{-0.15in} & .\hspace{-0.15in} & .\hspace{-0.15in} &
.\hspace{-0.15in} & .\hspace{-0.15in} & .\hspace{-0.15in} &
.\hspace{-0.15in}
& .\hspace{-0.15in} & .\hspace{-0.15in} & .\hspace{-0.15in} & 0 \\
\hspace{-0.15in} & \hspace{-0.15in} & \hspace{-0.15in} &
\hspace{-0.15in} & \hspace{-0.15in} & \hspace{-0.15in} &
\hspace{-0.15in} & \hspace{-0.15in} & \hspace{-0.15in} &
\hspace{-0.15in} & \hspace{-0.15in} & \hspace{-0.15in} &
\hspace{-0.15in} & \hspace{-0.15in} & \hspace{-0.15in} &
\hspace{-0.15in} & \hspace{-0.15in} & \hspace{-0.15in} &
\hspace{-0.15in} & \hspace{-0.15in} & \hspace{-0.15in} &
\hspace{-0.15in} & \hspace{-0.15in} & \hspace{-0.15in} &
\hspace{-0.15in} & \hspace{-0.15in} & \hspace{-0.15in} &
\hspace{-0.15in} &
\hspace{-0.15in} &  \\
\frac{1}{2}\hspace{-0.15in} & \frac{-1}{2}\hspace{-0.15in} & \frac{-1}{2}%
\hspace{-0.15in} & \frac{1}{2}\hspace{-0.15in} & \frac{1}{2}\hspace{-0.15in} & -%
\frac{1}{2}\hspace{-0.15in} & \frac{-1}{2}\hspace{-0.15in} & \frac{-1}{2}%
\hspace{-0.15in} & \frac{-1}{2}\hspace{-0.15in} &
\frac{1}{2}\hspace{-0.15in} &
\frac{1}{2}\hspace{-0.15in} & 1\hspace{-0.15in} & 0\hspace{-0.15in} & 0\hspace{%
-0.2in} & 0\hspace{-0.15in} & .\hspace{-0.15in} & .\hspace{-0.15in} & .\hspace{%
-0.2in} & .\hspace{-0.15in} & .\hspace{-0.15in} & .\hspace{-0.15in} & .\hspace{%
-0.2in} & .\hspace{-0.15in} & .\hspace{-0.15in} & .\hspace{-0.15in} & .\hspace{%
-0.2in} & .\hspace{-0.15in} & .\hspace{-0.15in} & .\hspace{-0.15in} & 0 \\
\hspace{-0.15in} & \hspace{-0.15in} & \hspace{-0.15in} &
\hspace{-0.15in} & \hspace{-0.15in} & \hspace{-0.15in} &
\hspace{-0.15in} & \hspace{-0.15in} & \hspace{-0.15in} &
\hspace{-0.15in} & \hspace{-0.15in} & \hspace{-0.15in} &
\hspace{-0.15in} & \hspace{-0.15in} & \hspace{-0.15in} &
\hspace{-0.15in} & \hspace{-0.15in} & \hspace{-0.15in} &
\hspace{-0.15in} & \hspace{-0.15in} & \hspace{-0.15in} &
\hspace{-0.15in} & \hspace{-0.15in} & \hspace{-0.15in} &
\hspace{-0.15in} & \hspace{-0.15in} & \hspace{-0.15in} &
\hspace{-0.15in} &
\hspace{-0.15in} &  \\
\frac{1}{12}\hspace{-0.15in} & \frac{-1}{12}\hspace{-0.15in} & \frac{-1}{12}%
\hspace{-0.15in} & \frac{-1}{12}\hspace{-0.15in} &
\frac{1}{12}\hspace{-0.15in}
& \frac{1}{12}\hspace{-0.15in} & \frac{-1}{12}\hspace{-0.15in} & \frac{3}{4}%
\hspace{-0.15in} & \frac{1}{4}\hspace{-0.15in} &
\frac{1}{4}\hspace{-0.15in} & \frac{1}{4}\hspace{-0.15in} &
0\hspace{-0.15in} & \frac{1}{2}\hspace{-0.15in} &
\frac{1}{2}\hspace{-0.15in} & \frac{-1}{2}\hspace{-0.15in} & \frac{-1}{2}%
\hspace{-0.15in} & -1\hspace{-0.15in} & 0\hspace{-0.15in} & 0\hspace{-0.15in} & 0%
\hspace{-0.15in} & .\hspace{-0.15in} & .\hspace{-0.15in} & .\hspace{-0.15in} & .%
\hspace{-0.15in} & .\hspace{-0.15in} & .\hspace{-0.15in} & .\hspace{-0.15in} & .%
\hspace{-0.15in} & .\hspace{-0.15in} & 0 \\
\hspace{-0.15in} & \hspace{-0.15in} & \hspace{-0.15in} &
\hspace{-0.15in} & \hspace{-0.15in} & \hspace{-0.15in} &
\hspace{-0.15in} & \hspace{-0.15in} & \hspace{-0.15in} &
\hspace{-0.15in} & \hspace{-0.15in} & \hspace{-0.15in} &
\hspace{-0.15in} & \hspace{-0.15in} & \hspace{-0.15in} &
\hspace{-0.15in} & \hspace{-0.15in} & \hspace{-0.15in} &
\hspace{-0.15in} & \hspace{-0.15in} & \hspace{-0.15in} &
\hspace{-0.15in} & \hspace{-0.15in} & \hspace{-0.15in} &
\hspace{-0.15in} & \hspace{-0.15in} & \hspace{-0.15in} &
\hspace{-0.15in} &
\hspace{-0.15in} &  \\
\frac{1}{12}\hspace{-0.15in} & \frac{-1}{12}\hspace{-0.15in} & \frac{-1}{12}%
\hspace{-0.15in} & \frac{-1}{12}\hspace{-0.15in} &
\frac{1}{12}\hspace{-0.15in}
& \frac{1}{12}\hspace{-0.15in} & \frac{-1}{12}\hspace{-0.15in} & \frac{-1}{4}%
\hspace{-0.15in} & \frac{-1}{4}\hspace{-0.15in} &
\frac{1}{4}\hspace{-0.15in} & \frac{1}{4}\hspace{-0.15in} &
0\hspace{-0.15in} & \frac{-1}{2}\hspace{-0.15in}
& \frac{-1}{2}\hspace{-0.15in} & \frac{1}{2}\hspace{-0.15in} & \frac{1}{2}%
\hspace{-0.15in} & 1\hspace{-0.15in} & 0\hspace{-0.15in} & 0\hspace{-0.15in} & 0%
\hspace{-0.15in} & .\hspace{-0.15in} & .\hspace{-0.15in} & .\hspace{-0.15in} & .%
\hspace{-0.15in} & .\hspace{-0.15in} & .\hspace{-0.15in} & .\hspace{-0.15in} & .%
\hspace{-0.15in} & .\hspace{-0.15in} & 0 \\
\hspace{-0.15in} & \hspace{-0.15in} & \hspace{-0.15in} &
\hspace{-0.15in} & \hspace{-0.15in} & \hspace{-0.15in} &
\hspace{-0.15in} & \hspace{-0.15in} & \hspace{-0.15in} &
\hspace{-0.15in} & \hspace{-0.15in} & \hspace{-0.15in} &
\hspace{-0.15in} & \hspace{-0.15in} & \hspace{-0.15in} &
\hspace{-0.15in} & \hspace{-0.15in} & \hspace{-0.15in} &
\hspace{-0.15in} & \hspace{-0.15in} & \hspace{-0.15in} &
\hspace{-0.15in} & \hspace{-0.15in} & \hspace{-0.15in} &
\hspace{-0.15in} & \hspace{-0.15in} & \hspace{-0.15in} &
\hspace{-0.15in} &
\hspace{-0.15in} &  \\
\frac{-1}{6}\hspace{-0.15in} & \frac{1}{6}\hspace{-0.15in} & \frac{1}{6}%
\hspace{-0.15in} & \frac{1}{6}\hspace{-0.15in} &
\frac{-1}{6}\hspace{-0.15in} &
\frac{-1}{6}\hspace{-0.15in} & \frac{1}{6}\hspace{-0.15in} & \frac{1}{2}%
\hspace{-0.15in} & \frac{1}{2}\hspace{-0.15in} &
\frac{-1}{2}\hspace{-0.15in} &
\frac{-1}{2}\hspace{-0.15in} & 0\hspace{-0.15in} & 0\hspace{-0.15in} & 0\hspace{%
-0.2in} & 0\hspace{-0.15in} & 0\hspace{-0.15in} & 0\hspace{-0.15in} & 0\hspace{%
-0.2in} & 0\hspace{-0.15in} & 0\hspace{-0.15in} & \hspace{-0.15in} & \hspace{%
-0.2in} & \hspace{-0.15in} & \hspace{-0.15in} & \hspace{-0.15in} & \hspace{%
-0.2in} & \hspace{-0.15in} & \hspace{-0.15in} & \hspace{-0.15in} &  \\
\hspace{-0.15in} & \hspace{-0.15in} & \hspace{-0.15in} &
\hspace{-0.15in} & \hspace{-0.15in} & \hspace{-0.15in} &
\hspace{-0.15in} & \hspace{-0.15in} & \hspace{-0.15in} &
\hspace{-0.15in} & \hspace{-0.15in} & \hspace{-0.15in} &
\hspace{-0.15in} & \hspace{-0.15in} & \hspace{-0.15in} &
\hspace{-0.15in} & \hspace{-0.15in} & \hspace{-0.15in} &
\hspace{-0.15in} & \hspace{-0.15in} & \hspace{-0.15in} &
\hspace{-0.15in} & \hspace{-0.15in} & \hspace{-0.15in} &
\hspace{-0.15in} & \hspace{-0.15in} & \hspace{-0.15in} &
\hspace{-0.15in} &
\hspace{-0.15in} &  \\
\frac{-1}{12}\hspace{-0.15in} & \frac{1}{12}\hspace{-0.15in} & \frac{1}{12}%
\hspace{-0.15in} & \frac{1}{12}\hspace{-0.15in} &
\frac{-1}{12}\hspace{-0.15in}
& \frac{-1}{12}\hspace{-0.15in} & \frac{1}{12}\hspace{-0.15in} & \frac{1}{4}%
\hspace{-0.15in} & \frac{1}{4}\hspace{-0.15in} &
\frac{-1}{4}\hspace{-0.15in} & \frac{3}{4}\hspace{-0.15in} &
0\hspace{-0.15in} & \frac{1}{2}\hspace{-0.15in} &
\frac{1}{2}\hspace{-0.15in} & \frac{-1}{2}\hspace{-0.15in} & \frac{-1}{2}%
\hspace{-0.15in} & -1\hspace{-0.15in} & 0\hspace{-0.15in} & 0\hspace{-0.15in} & 0%
\hspace{-0.15in} & .\hspace{-0.15in} & .\hspace{-0.15in} & .\hspace{-0.15in} & .%
\hspace{-0.15in} & .\hspace{-0.15in} & .\hspace{-0.15in} & .\hspace{-0.15in} & .%
\hspace{-0.15in} & .\hspace{-0.15in} & 0 \\
\hspace{-0.15in} & \hspace{-0.15in} & \hspace{-0.15in} &
\hspace{-0.15in} & \hspace{-0.15in} & \hspace{-0.15in} &
\hspace{-0.15in} & \hspace{-0.15in} & \hspace{-0.15in} &
\hspace{-0.15in} & \hspace{-0.15in} & \hspace{-0.15in} &
\hspace{-0.15in} & \hspace{-0.15in} & \hspace{-0.15in} &
\hspace{-0.15in} & \hspace{-0.15in} & \hspace{-0.15in} &
\hspace{-0.15in} & \hspace{-0.15in} & \hspace{-0.15in} &
\hspace{-0.15in} & \hspace{-0.15in} & \hspace{-0.15in} &
\hspace{-0.15in} & \hspace{-0.15in} & \hspace{-0.15in} &
\hspace{-0.15in} &
\hspace{-0.15in} &  \\
\frac{1}{4}\hspace{-0.15in} & \frac{-1}{4}\hspace{-0.15in} & \frac{-1}{4}%
\hspace{-0.15in} & \frac{-1}{4}\hspace{-0.15in} &
\frac{1}{4}\hspace{-0.15in} &
\frac{1}{4}\hspace{-0.15in} & \frac{-1}{4}\hspace{-0.15in} & \frac{-3}{4}%
\hspace{-0.15in} & \frac{1}{4}\hspace{-0.15in} & \frac{3}{4}\hspace{-0.15in} & -%
\frac{1}{4}\hspace{-0.15in} & 0\hspace{-0.15in} &
\frac{-1}{2}\hspace{-0.15in}
& \frac{-1}{2}\hspace{-0.15in} & \frac{1}{2}\hspace{-0.15in} & \frac{1}{2}%
\hspace{-0.15in} & 1\hspace{-0.15in} & 0\hspace{-0.15in} & 0\hspace{-0.15in} & 0%
\hspace{-0.15in} & .\hspace{-0.15in} & .\hspace{-0.15in} & .\hspace{-0.15in} & .%
\hspace{-0.15in} & .\hspace{-0.15in} & .\hspace{-0.15in} & .\hspace{-0.15in} & .%
\hspace{-0.15in} & .\hspace{-0.15in} & 0 \\
\hspace{-0.15in} & \hspace{-0.15in} & \hspace{-0.15in} &
\hspace{-0.15in} & \hspace{-0.15in} & \hspace{-0.15in} &
\hspace{-0.15in} & \hspace{-0.15in} & \hspace{-0.15in} &
\hspace{-0.15in} & \hspace{-0.15in} & \hspace{-0.15in} &
\hspace{-0.15in} & \hspace{-0.15in} & \hspace{-0.15in} &
\hspace{-0.15in} & \hspace{-0.15in} & \hspace{-0.15in} &
\hspace{-0.15in} & \hspace{-0.15in} & \hspace{-0.15in} &
\hspace{-0.15in} & \hspace{-0.15in} & \hspace{-0.15in} &
\hspace{-0.15in} & \hspace{-0.15in} & \hspace{-0.15in} &
\hspace{-0.15in} &
\hspace{-0.15in} &  \\
\frac{1}{24}\hspace{-0.15in} & \frac{-1}{24}\hspace{-0.15in} & \frac{-1}{24}%
\hspace{-0.15in} & \frac{-1}{24}\hspace{-0.15in} &
\frac{1}{24}\hspace{-0.15in}
& \frac{1}{24}\hspace{-0.15in} & \frac{-1}{24}\hspace{-0.15in} & \frac{-1}{8}%
\hspace{-0.15in} & \frac{-1}{8}\hspace{-0.15in} &
\frac{1}{8}\hspace{-0.15in} & \frac{1}{8}\hspace{-0.15in} &
0\hspace{-0.15in} & \frac{3}{4}\hspace{-0.15in} &
\frac{1}{4}\hspace{-0.15in} & \frac{1}{4}\hspace{-0.15in} & \frac{1}{4}\hspace{%
-0.2in} & 0\hspace{-0.15in} & \frac{1}{2}\hspace{-0.15in} & \frac{1}{2}\hspace{%
-0.2in} & \frac{-1}{2}\hspace{-0.15in} & \frac{-1}{2}\hspace{-0.15in} & -1%
\hspace{-0.15in} & 0\hspace{-0.15in} & .\hspace{-0.15in} & .\hspace{-0.15in} & .%
\hspace{-0.15in} & .\hspace{-0.15in} & .\hspace{-0.15in} &
.\hspace{-0.15in} & .
\\
\hspace{-0.15in} & \hspace{-0.15in} & \hspace{-0.15in} &
\hspace{-0.15in} & \hspace{-0.15in} & \hspace{-0.15in} &
\hspace{-0.15in} & \hspace{-0.15in} & \hspace{-0.15in} &
\hspace{-0.15in} & \hspace{-0.15in} & \hspace{-0.15in} &
\hspace{-0.15in} & \hspace{-0.15in} & \hspace{-0.15in} &
\hspace{-0.15in} & \hspace{-0.15in} & \hspace{-0.15in} &
\hspace{-0.15in} & \hspace{-0.15in} & \hspace{-0.15in} &
\hspace{-0.15in} & \hspace{-0.15in} & \hspace{-0.15in} &
\hspace{-0.15in} & \hspace{-0.15in} & \hspace{-0.15in} &
\hspace{-0.15in} &
\hspace{-0.15in} &  \\
\frac{1}{24}\hspace{-0.15in} & \frac{-1}{24}\hspace{-0.15in} & \frac{-1}{24}%
\hspace{-0.15in} & \frac{-1}{24}\hspace{-0.15in} &
\frac{1}{24}\hspace{-0.15in}
& \frac{1}{24}\hspace{-0.15in} & \frac{-1}{24}\hspace{-0.15in} & \frac{-1}{8}%
\hspace{-0.15in} & \frac{-1}{8}\hspace{-0.15in} &
\frac{1}{8}\hspace{-0.15in} & \frac{1}{8}\hspace{-0.15in} &
0\hspace{-0.15in} & \frac{-1}{4}\hspace{-0.15in}
& \frac{-1}{4}\hspace{-0.15in} & \frac{1}{4}\hspace{-0.15in} & \frac{1}{4}%
\hspace{-0.15in} & 0\hspace{-0.15in} & \frac{-1}{2}\hspace{-0.15in} & \frac{-1}{%
2}\hspace{-0.15in} & \frac{1}{2}\hspace{-0.15in} &
\frac{1}{2}\hspace{-0.15in} & 1\hspace{-0.15in} &
0\hspace{-0.15in} & .\hspace{-0.15in} & .\hspace{-0.15in} &
.\hspace{-0.15in} & .\hspace{-0.15in} & .\hspace{-0.15in} &
.\hspace{-0.15in}
& 0 \\
\hspace{-0.15in} & \hspace{-0.15in} & \hspace{-0.15in} &
\hspace{-0.15in} & \hspace{-0.15in} & \hspace{-0.15in} &
\hspace{-0.15in} & \hspace{-0.15in} & \hspace{-0.15in} &
\hspace{-0.15in} & \hspace{-0.15in} & \hspace{-0.15in} &
\hspace{-0.15in} & \hspace{-0.15in} & \hspace{-0.15in} &
\hspace{-0.15in} & \hspace{-0.15in} & \hspace{-0.15in} &
\hspace{-0.15in} & \hspace{-0.15in} & \hspace{-0.15in} &
\hspace{-0.15in} & \hspace{-0.15in} & \hspace{-0.15in} &
\hspace{-0.15in} & \hspace{-0.15in} & \hspace{-0.15in} &
\hspace{-0.15in} &
\hspace{-0.15in} &  \\
\frac{-1}{12}\hspace{-0.15in} & \frac{1}{12}\hspace{-0.15in} & \frac{1}{12}%
\hspace{-0.15in} & \frac{1}{12}\hspace{-0.15in} &
\frac{-1}{12}\hspace{-0.15in}
& \frac{-1}{12}\hspace{-0.15in} & \frac{1}{12}\hspace{-0.15in} & \frac{1}{4}%
\hspace{-0.15in} & \frac{1}{4}\hspace{-0.15in} &
\frac{-1}{4}\hspace{-0.15in} & \frac{-1}{4}\hspace{-0.15in} &
0\hspace{-0.15in} & \frac{1}{2}\hspace{-0.15in}
& \frac{1}{2}\hspace{-0.15in} & \frac{-1}{2}\hspace{-0.15in} & \frac{-1}{2}%
\hspace{-0.15in} & 0\hspace{-0.15in} & 0\hspace{-0.15in} & 0\hspace{-0.15in} & 0%
\hspace{-0.15in} & 0\hspace{-0.15in} & 0\hspace{-0.15in} & 0\hspace{-0.15in} & .%
\hspace{-0.15in} & .\hspace{-0.15in} & .\hspace{-0.15in} & .\hspace{-0.15in} & .%
\hspace{-0.15in} & .\hspace{-0.15in} & 0 \\
\hspace{-0.15in} & \hspace{-0.15in} & \hspace{-0.15in} &
\hspace{-0.15in} & \hspace{-0.15in} & \hspace{-0.15in} &
\hspace{-0.15in} & \hspace{-0.15in} & \hspace{-0.15in} &
\hspace{-0.15in} & \hspace{-0.15in} & \hspace{-0.15in} &
\hspace{-0.15in} & \hspace{-0.15in} & \hspace{-0.15in} &
\hspace{-0.15in} & \hspace{-0.15in} & \hspace{-0.15in} &
\hspace{-0.15in} & \hspace{-0.15in} & \hspace{-0.15in} &
\hspace{-0.15in} & \hspace{-0.15in} & \hspace{-0.15in} &
\hspace{-0.15in} & \hspace{-0.15in} & \hspace{-0.15in} &
\hspace{-0.15in} &
\hspace{-0.15in} &  \\
\frac{-1}{24}\hspace{-0.15in} & \frac{1}{24}\hspace{-0.15in} & \frac{1}{24}%
\hspace{-0.15in} & \frac{1}{24}\hspace{-0.15in} &
\frac{-1}{24}\hspace{-0.15in}
& \frac{-1}{24}\hspace{-0.15in} & \frac{1}{24}\hspace{-0.15in} & \frac{1}{8}%
\hspace{-0.15in} & \frac{1}{8}\hspace{-0.15in} &
\frac{-1}{8}\hspace{-0.15in} & \frac{-1}{8}\hspace{-0.15in} &
0\hspace{-0.15in} & \frac{1}{4}\hspace{-0.15in}
& \frac{1}{4}\hspace{-0.15in} & \frac{-1}{4}\hspace{-0.15in} & \frac{3}{4}%
\hspace{-0.15in} & 0\hspace{-0.15in} & \frac{1}{2}\hspace{-0.15in} & \frac{1}{2}%
\hspace{-0.15in} & \frac{-1}{2}\hspace{-0.15in} &
\frac{-1}{2}\hspace{-0.15in} & -1\hspace{-0.15in} &
0\hspace{-0.15in} & .\hspace{-0.15in} & .\hspace{-0.15in} &
.\hspace{-0.15in} & .\hspace{-0.15in} & .\hspace{-0.15in} &
.\hspace{-0.15in}
& 0 \\
\hspace{-0.15in} & \hspace{-0.15in} & \hspace{-0.15in} &
\hspace{-0.15in} & \hspace{-0.15in} & \hspace{-0.15in} &
\hspace{-0.15in} & \hspace{-0.15in} & \hspace{-0.15in} &
\hspace{-0.15in} & \hspace{-0.15in} & \hspace{-0.15in} &
\hspace{-0.15in} & \hspace{-0.15in} & \hspace{-0.15in} &
\hspace{-0.15in} & \hspace{-0.15in} & \hspace{-0.15in} &
\hspace{-0.15in} & \hspace{-0.15in} & \hspace{-0.15in} &
\hspace{-0.15in} & \hspace{-0.15in} & \hspace{-0.15in} &
\hspace{-0.15in} & \hspace{-0.15in} & \hspace{-0.15in} &
\hspace{-0.15in} &
\hspace{-0.15in} &  \\
\frac{1}{8}\hspace{-0.15in} & \frac{-1}{8}\hspace{-0.15in} & \frac{-1}{8}%
\hspace{-0.15in} & \frac{-1}{8}\hspace{-0.15in} &
\frac{1}{8}\hspace{-0.15in} &
\frac{1}{8}\hspace{-0.15in} & \frac{-1}{8}\hspace{-0.15in} & \frac{-3}{8}%
\hspace{-0.15in} & \frac{-3}{8}\hspace{-0.15in} &
\frac{3}{8}\hspace{-0.15in} & \frac{3}{8}\hspace{-0.15in} &
0\hspace{-0.15in} & \frac{-3}{4}\hspace{-0.15in}
& \frac{1}{4}\hspace{-0.15in} & \frac{3}{4}\hspace{-0.15in} & \frac{-1}{4}%
\hspace{-0.15in} & 0\hspace{-0.15in} & \frac{-1}{2}\hspace{-0.15in} & \frac{-1}{%
2}\hspace{-0.15in} & \frac{1}{2}\hspace{-0.15in} &
\frac{1}{2}\hspace{-0.15in} & 1\hspace{-0.15in} &
0\hspace{-0.15in} & .\hspace{-0.15in} & .\hspace{-0.15in} &
.\hspace{-0.15in} & .\hspace{-0.15in} & .\hspace{-0.15in} &
.\hspace{-0.15in}
& 0 \\
\hspace{-0.15in} & \hspace{-0.15in} & \hspace{-0.15in} &
\hspace{-0.15in} & \hspace{-0.15in} & \hspace{-0.15in} &
\hspace{-0.15in} & \hspace{-0.15in} & \hspace{-0.15in} &
\hspace{-0.15in} & \hspace{-0.15in} & \hspace{-0.15in} &
\hspace{-0.15in} & \hspace{-0.15in} & \hspace{-0.15in} &
\hspace{-0.15in} & \hspace{-0.15in} & \hspace{-0.15in} &
\hspace{-0.15in} & \hspace{-0.15in} & \hspace{-0.15in} &
\hspace{-0.15in} & \hspace{-0.15in} & \hspace{-0.15in} &
\hspace{-0.15in} & \hspace{-0.15in} & \hspace{-0.15in} &
\hspace{-0.15in} &
\hspace{-0.15in} &  \\
.\hspace{-0.15in} & .\hspace{-0.15in} & .\hspace{-0.15in} & .\hspace{-0.15in} & .%
\hspace{-0.15in} & .\hspace{-0.15in} & .\hspace{-0.15in} & .\hspace{-0.15in} & .%
\hspace{-0.15in} & .\hspace{-0.15in} & .\hspace{-0.15in} & .\hspace{-0.15in} & .%
\hspace{-0.15in} & .\hspace{-0.15in} & .\hspace{-0.15in} & ..\hspace{-0.15in} & .%
\hspace{-0.15in} & .\hspace{-0.15in} & .\hspace{-0.15in} & .\hspace{-0.15in} & .%
\hspace{-0.15in} & .\hspace{-0.15in} & .\hspace{-0.15in} & .\hspace{-0.15in} & .%
\hspace{-0.15in} & .\hspace{-0.15in} & .\hspace{-0.15in} & .\hspace{-0.15in} & .%
\hspace{-0.15in} & 0 \\
\hspace{-0.15in} & \hspace{-0.15in} & \hspace{-0.15in} &
\hspace{-0.15in} & \hspace{-0.15in} & \hspace{-0.15in} &
\hspace{-0.15in} & \hspace{-0.15in} & \hspace{-0.15in} &
\hspace{-0.15in} & \hspace{-0.15in} & \hspace{-0.15in} &
\hspace{-0.15in} & \hspace{-0.15in} & \hspace{-0.15in} &
\hspace{-0.15in} & \hspace{-0.15in} & \hspace{-0.15in} &
\hspace{-0.15in} & \hspace{-0.15in} & \hspace{-0.15in} &
\hspace{-0.15in} & \hspace{-0.15in} & \hspace{-0.15in} &
\hspace{-0.15in} & \hspace{-0.15in} & \hspace{-0.15in} &
\hspace{-0.15in} &
\hspace{-0.15in} &  \\
.\hspace{-0.15in} & .\hspace{-0.15in} & .\hspace{-0.15in} & .\hspace{-0.15in} & .%
\hspace{-0.15in} & .\hspace{-0.15in} & .\hspace{-0.15in} & .\hspace{-0.15in} & .%
\hspace{-0.15in} & .\hspace{-0.15in} & .\hspace{-0.15in} & .\hspace{-0.15in} & .%
\hspace{-0.15in} & .\hspace{-0.15in} & .\hspace{-0.15in} & .\hspace{-0.15in} & .%
\hspace{-0.15in} & .\hspace{-0.15in} & .\hspace{-0.15in} & .\hspace{-0.15in} & .%
\hspace{-0.15in} & .\hspace{-0.15in} & .\hspace{-0.15in} & .\hspace{-0.15in} & .%
\hspace{-0.15in} & .\hspace{-0.15in} & .\hspace{-0.15in} & .\hspace{-0.15in} & .%
\hspace{-0.15in} & 0 \\
\hspace{-0.15in} & \hspace{-0.15in} & \hspace{-0.15in} &
\hspace{-0.15in} & \hspace{-0.15in} & \hspace{-0.15in} &
\hspace{-0.15in} & \hspace{-0.15in} & \hspace{-0.15in} &
\hspace{-0.15in} & \hspace{-0.15in} & \hspace{-0.15in} &
\hspace{-0.15in} & \hspace{-0.15in} & \hspace{-0.15in} &
\hspace{-0.15in} & \hspace{-0.15in} & \hspace{-0.15in} &
\hspace{-0.15in} & \hspace{-0.15in} & \hspace{-0.15in} &
\hspace{-0.15in} & \hspace{-0.15in} & \hspace{-0.15in} &
\hspace{-0.15in} & \hspace{-0.15in} & \hspace{-0.15in} &
\hspace{-0.15in} &
\hspace{-0.15in} &  \\
\frac{1}{2^{n-1}\times 3}\hspace{-0.15in} & \frac{-1}{2^{n-1}\times 3}\hspace{%
-0.2in} & \frac{-1}{2^{n-1}\times 3}\hspace{-0.15in} & \frac{-1}{%
2^{n-1}\times 3}\hspace{-0.15in} & \frac{1}{2^{n-1}\times
3}\hspace{-0.15in} &
\frac{1}{2^{n-1}\times 3}\hspace{-0.15in} & \frac{-1}{2^{n-1}\times 3}\hspace{%
-0.2in} & \frac{-1}{2^{n}}\hspace{-0.15in} &
\frac{-1}{2^{n}}\hspace{-0.15in}
& \frac{1}{2^{n}}\hspace{-0.15in} & \frac{1}{2^{n}}\hspace{-0.15in} & 0\hspace{%
-0.2in} & .\hspace{-0.15in} & .\hspace{-0.15in} &
\frac{-1}{8}\hspace{-0.15in}
& \frac{-1}{8}\hspace{-0.15in} & \frac{1}{8}\hspace{-0.15in} & \frac{1}{8}%
\hspace{-0.15in} & 0\hspace{-0.15in} & \frac{3}{4}\hspace{-0.15in} & \frac{1}{4}%
\hspace{-0.15in} & \frac{1}{4}\hspace{-0.15in} & \frac{1}{4}\hspace{-0.15in} & 0%
\hspace{-0.15in} & \frac{1}{2}\hspace{-0.15in} & \frac{1}{2}\hspace{-0.15in} & %
\frac{-1}{2}\hspace{-0.15in} & \frac{-1}{2}\hspace{-0.15in} &
-1\hspace{-0.15in}
& 0 \\
\hspace{-0.15in} & \hspace{-0.15in} & \hspace{-0.15in} &
\hspace{-0.15in} & \hspace{-0.15in} & \hspace{-0.15in} &
\hspace{-0.15in} & \hspace{-0.15in} & \hspace{-0.15in} &
\hspace{-0.15in} & \hspace{-0.15in} & \hspace{-0.15in} &
\hspace{-0.15in} & \hspace{-0.15in} & \hspace{-0.15in} &
\hspace{-0.15in} & \hspace{-0.15in} & \hspace{-0.15in} &
\hspace{-0.15in} & \hspace{-0.15in} & \hspace{-0.15in} &
\hspace{-0.15in} & \hspace{-0.15in} & \hspace{-0.15in} &
\hspace{-0.15in} & \hspace{-0.15in} & \hspace{-0.15in} &
\hspace{-0.15in} &
\hspace{-0.15in} &  \\
\frac{1}{2^{n-1}\times 3}\hspace{-0.15in} & \frac{-1}{2^{n-1}\times 3}\hspace{%
-0.2in} & \frac{-1}{2^{n-1}\times 3}\hspace{-0.15in} & \frac{-1}{%
2^{n-1}\times 3}\hspace{-0.15in} & \frac{1}{2^{n-1}\times
3}\hspace{-0.15in} &
\frac{1}{2^{n-1}\times 3}\hspace{-0.15in} & \frac{-1}{2^{n-1}\times 3}\hspace{%
-0.2in} & \frac{-1}{2^{n}}\hspace{-0.15in} &
\frac{-1}{2^{n}}\hspace{-0.15in}
& \frac{1}{2^{n}}\hspace{-0.15in} & \frac{1}{2^{n}}\hspace{-0.15in} & 0\hspace{%
-0.2in} & .\hspace{-0.15in} & .\hspace{-0.15in} &
\frac{-1}{8}\hspace{-0.15in}
& \frac{-1}{8}\hspace{-0.15in} & \frac{1}{8}\hspace{-0.15in} & \frac{1}{8}%
\hspace{-0.15in} & 0\hspace{-0.15in} & \frac{-1}{4}\hspace{-0.15in} & \frac{-1}{%
4}\hspace{-0.15in} & \frac{1}{4}\hspace{-0.15in} &
\frac{1}{4}\hspace{-0.15in}
& 0\hspace{-0.15in} & \frac{-1}{2}\hspace{-0.15in} & \frac{-1}{2}\hspace{-0.2in%
} & \frac{1}{2}\hspace{-0.15in} & \frac{1}{2}\hspace{-0.15in} & 1\hspace{-0.2in%
} & 0 \\
\hspace{-0.15in} & \hspace{-0.15in} & \hspace{-0.15in} &
\hspace{-0.15in} & \hspace{-0.15in} & \hspace{-0.15in} &
\hspace{-0.15in} & \hspace{-0.15in} & \hspace{-0.15in} &
\hspace{-0.15in} & \hspace{-0.15in} & \hspace{-0.15in} &
\hspace{-0.15in} & \hspace{-0.15in} & \hspace{-0.15in} &
\hspace{-0.15in} & \hspace{-0.15in} & \hspace{-0.15in} &
\hspace{-0.15in} & \hspace{-0.15in} & \hspace{-0.15in} &
\hspace{-0.15in} & \hspace{-0.15in} & \hspace{-0.15in} &
\hspace{-0.15in} & \hspace{-0.15in} & \hspace{-0.15in} &
\hspace{-0.15in} &
\hspace{-0.15in} &  \\
\frac{-1}{2^{n-1}\times 3}\hspace{-0.15in} & \frac{1}{2^{n-1}\times 3}\hspace{%
-0.2in} & \frac{1}{2^{n-1}\times 3}\hspace{-0.15in} & \frac{1}{2^{n-1}\times 3%
}\hspace{-0.15in} & \frac{-1}{2^{n-1}\times 3}\hspace{-0.15in} & \frac{-1}{%
2^{n-1}\times 3}\hspace{-0.15in} & \frac{1}{2^{n-1}\times
3}\hspace{-0.15in} &
\frac{1}{2^{n-1}}\hspace{-0.15in} & \frac{1}{2^{n-1}}\hspace{-0.15in} & \frac{-%
1}{2^{n-1}}\hspace{-0.15in} & \frac{-1}{2^{n-1}}\hspace{-0.15in} & 0\hspace{%
-0.2in} & .\hspace{-0.15in} & .\hspace{-0.15in} &
\frac{1}{4}\hspace{-0.15in} &
\frac{1}{4}\hspace{-0.15in} & \frac{-1}{4}\hspace{-0.15in} & \frac{-1}{4}%
\hspace{-0.15in} & 0\hspace{-0.15in} & \frac{1}{2}\hspace{-0.15in} & \frac{1}{2}%
\hspace{-0.15in} & \frac{-1}{2}\hspace{-0.15in} &
\frac{-1}{2}\hspace{-0.15in} & 0\hspace{-0.15in} &
0\hspace{-0.15in} & 0\hspace{-0.15in} & 0\hspace{-0.15in}
& 0\hspace{-0.15in} & 0\hspace{-0.15in} & 0 \\
\hspace{-0.15in} & \hspace{-0.15in} & \hspace{-0.15in} &
\hspace{-0.15in} & \hspace{-0.15in} & \hspace{-0.15in} &
\hspace{-0.15in} & \hspace{-0.15in} & \hspace{-0.15in} &
\hspace{-0.15in} & \hspace{-0.15in} & \hspace{-0.15in} &
\hspace{-0.15in} & \hspace{-0.15in} & \hspace{-0.15in} &
\hspace{-0.15in} & \hspace{-0.15in} & \hspace{-0.15in} &
\hspace{-0.15in} & \hspace{-0.15in} & \hspace{-0.15in} &
\hspace{-0.15in} & \hspace{-0.15in} & \hspace{-0.15in} &
\hspace{-0.15in} & \hspace{-0.15in} & \hspace{-0.15in} &
\hspace{-0.15in} &
\hspace{-0.15in} &  \\
\frac{-1}{2^{n-1}\times 3}\hspace{-0.15in} & \frac{1}{2^{n-1}\times 3}\hspace{%
-0.2in} & \frac{1}{2^{n-1}\times 3}\hspace{-0.15in} & \frac{1}{2^{n-1}\times 3%
}\hspace{-0.15in} & \frac{-1}{2^{n-1}\times 3}\hspace{-0.15in} & \frac{-1}{%
2^{n-1}\times 3}\hspace{-0.15in} & \frac{1}{2^{n-1}\times
3}\hspace{-0.15in} &
\frac{1}{2^{n}}\hspace{-0.15in} & \frac{1}{2^{n}}\hspace{-0.15in} & \frac{-1}{%
2^{n}}\hspace{-0.15in} & \frac{-1}{2^{n}}\hspace{-0.15in} &
0\hspace{-0.15in} &
.\hspace{-0.15in} & .\hspace{-0.15in} & \frac{1}{8}\hspace{-0.15in} & \frac{1}{8%
}\hspace{-0.15in} & \frac{-1}{8}\hspace{-0.15in} &
\frac{-1}{8}\hspace{-0.15in} & 0\hspace{-0.15in} &
\frac{1}{4}\hspace{-0.15in} & \frac{1}{4}\hspace{-0.15in} &
\frac{-1}{4}\hspace{-0.15in} & \frac{3}{4}\hspace{-0.15in} &
0\hspace{-0.15in}
& \frac{1}{2}\hspace{-0.15in} & \frac{1}{2}\hspace{-0.15in} & \frac{-1}{2}%
\hspace{-0.15in} & \frac{-1}{2}\hspace{-0.15in} & -1\hspace{-0.15in} & 0 \\
\hspace{-0.15in} & \hspace{-0.15in} & \hspace{-0.15in} &
\hspace{-0.15in} & \hspace{-0.15in} & \hspace{-0.15in} &
\hspace{-0.15in} & \hspace{-0.15in} & \hspace{-0.15in} &
\hspace{-0.15in} & \hspace{-0.15in} & \hspace{-0.15in} &
\hspace{-0.15in} & \hspace{-0.15in} & \hspace{-0.15in} &
\hspace{-0.15in} & \hspace{-0.15in} & \hspace{-0.15in} &
\hspace{-0.15in} & \hspace{-0.15in} & \hspace{-0.15in} &
\hspace{-0.15in} & \hspace{-0.15in} & \hspace{-0.15in} &
\hspace{-0.15in} & \hspace{-0.15in} & \hspace{-0.15in} &
\hspace{-0.15in} &
\hspace{-0.15in} &  \\
\frac{1}{2^{n}}\hspace{-0.15in} & \frac{-1}{2^{n}}\hspace{-0.15in} & \frac{-1}{%
2^{n}}\hspace{-0.15in} & \frac{-1}{2^{n}}\hspace{-0.15in} & \frac{1}{2^{n}}%
\hspace{-0.15in} & \frac{1}{2^{n}}\hspace{-0.15in} & \frac{-1}{2^{n}}\hspace{%
-0.2in} & \frac{-3}{2^{n}}\hspace{-0.15in} &
\frac{-3}{2^{n}}\hspace{-0.15in}
& \frac{3}{2^{n}}\hspace{-0.15in} & \frac{3}{2^{n}}\hspace{-0.15in} & 0\hspace{%
-0.2in} & .\hspace{-0.15in} & .\hspace{-0.15in} &
\frac{-3}{8}\hspace{-0.15in}
& \frac{-3}{8}\hspace{-0.15in} & \frac{3}{8}\hspace{-0.15in} & \frac{3}{8}%
\hspace{-0.15in} & 0\hspace{-0.15in} & \frac{-3}{4}\hspace{-0.15in} & \frac{1}{4%
}\hspace{-0.15in} & \frac{3}{4}\hspace{-0.15in} &
\frac{-1}{4}\hspace{-0.15in}
& 0\hspace{-0.15in} & \frac{-1}{2}\hspace{-0.15in} & \frac{-1}{2}\hspace{-0.2in%
} & \frac{1}{2}\hspace{-0.15in} & \frac{1}{2}\hspace{-0.15in} & 1\hspace{-0.2in%
} & 0 \\
\hspace{-0.15in} & \hspace{-0.15in} & \hspace{-0.15in} &
\hspace{-0.15in} & \hspace{-0.15in} & \hspace{-0.15in} &
\hspace{-0.15in} & \hspace{-0.15in} & \hspace{-0.15in} &
\hspace{-0.15in} & \hspace{-0.15in} & \hspace{-0.15in} &
\hspace{-0.15in} & \hspace{-0.15in} & \hspace{-0.15in} &
\hspace{-0.15in} & \hspace{-0.15in} & \hspace{-0.15in} &
\hspace{-0.15in} & \hspace{-0.15in} & \hspace{-0.15in} &
\hspace{-0.15in} & \hspace{-0.15in} & \hspace{-0.15in} &
\hspace{-0.15in} & \hspace{-0.15in} & \hspace{-0.15in} &
\hspace{-0.15in} &
\hspace{-0.15in} &  \\
\frac{2^{n-1}+1}{2^{n-1}\times 3}\hspace{-0.15in} & \frac{-(2^{n-1}+1)}{%
2^{n-1}\times 3}\hspace{-0.15in} & \frac{2^{n}-1}{2^{n-1}\times 3}\hspace{%
-0.2in} & \frac{2^{n}-1}{2^{n-1}\times 3}\hspace{-0.15in} & \frac{-(2^{n}-1)}{%
2^{n-1}\times 3}\hspace{-0.15in} & \frac{-(2^{n}-1)}{2^{n-1}\times 3}\hspace{%
-0.2in} & \frac{-(2^{n+1}+1)}{2^{n-1}\times 3}\hspace{-0.15in} & \frac{-1}{2^{n}%
}\hspace{-0.15in} & \frac{-1}{2^{n}}\hspace{-0.15in} & \frac{1}{2^{n}}\hspace{%
-0.2in} & \frac{1}{2^{n}}\hspace{-0.15in} & 0\hspace{-0.15in} & .\hspace{-0.2in%
} & .\hspace{-0.15in} & \frac{-1}{8}\hspace{-0.15in} & \frac{-1}{8}\hspace{%
-0.2in} & \frac{1}{8}\hspace{-0.15in} & \frac{1}{8}\hspace{-0.15in} & 0\hspace{%
-0.2in} & \frac{-1}{4}\hspace{-0.15in} & \frac{-1}{4}\hspace{-0.15in} & \frac{1%
}{4}\hspace{-0.15in} & \frac{1}{4}\hspace{-0.15in} & 0\hspace{-0.15in} & \frac{1%
}{2}\hspace{-0.15in} & \frac{-1}{2}\hspace{-0.15in} & \frac{1}{2}\hspace{-0.2in%
} & \frac{-1}{2}\hspace{-0.15in} & 0\hspace{-0.15in} & 1 \\
\hspace{-0.15in} & \hspace{-0.15in} & \hspace{-0.15in} &
\hspace{-0.15in} & \hspace{-0.15in} & \hspace{-0.15in} &
\hspace{-0.15in} & \hspace{-0.15in} & \hspace{-0.15in} &
\hspace{-0.15in} & \hspace{-0.15in} & \hspace{-0.15in} &
\hspace{-0.15in} & \hspace{-0.15in} & \hspace{-0.15in} &
\hspace{-0.15in} & \hspace{-0.15in} & \hspace{-0.15in} &
\hspace{-0.15in} & \hspace{-0.15in} & \hspace{-0.15in} &
\hspace{-0.15in} & \hspace{-0.15in} & \hspace{-0.15in} &
\hspace{-0.15in} & \hspace{-0.15in} & \hspace{-0.15in} &
\hspace{-0.15in} &
\hspace{-0.15in} &  \\
\frac{-1}{3}\hspace{-0.15in} & \frac{1}{3}\hspace{-0.15in} & \frac{-2}{3}%
\hspace{-0.15in} & \frac{-2}{3}\hspace{-0.15in} &
\frac{2}{3}\hspace{-0.15in} & \frac{2}{3}\hspace{-0.15in} &
\frac{4}{3}\hspace{-0.15in} & 0\hspace{-0.15in} &
0\hspace{-0.15in} & 0\hspace{-0.15in} & 0\hspace{-0.15in} & 0\hspace{-0.15in} & .%
\hspace{-0.15in} & .\hspace{-0.15in} & 0\hspace{-0.15in} & 0\hspace{-0.15in} & 0%
\hspace{-0.15in} & 0\hspace{-0.15in} & 0\hspace{-0.15in} & 0\hspace{-0.15in} & 0%
\hspace{-0.15in} & 0\hspace{-0.15in} & 0\hspace{-0.15in} & 0\hspace{-0.15in} & 0%
\hspace{-0.15in} & 0\hspace{-0.15in} & 0\hspace{-0.15in} & 1\hspace{-0.15in} & 0%
\hspace{-0.15in} & -1 \\
\hspace{-0.15in} & \hspace{-0.15in} & \hspace{-0.15in} &
\hspace{-0.15in} & \hspace{-0.15in} & \hspace{-0.15in} &
\hspace{-0.15in} & \hspace{-0.15in} & \hspace{-0.15in} &
\hspace{-0.15in} & \hspace{-0.15in} & \hspace{-0.15in} &
\hspace{-0.15in} & \hspace{-0.15in} & \hspace{-0.15in} &
\hspace{-0.15in} & \hspace{-0.15in} & \hspace{-0.15in} &
\hspace{-0.15in} & \hspace{-0.15in} & \hspace{-0.15in} &
\hspace{-0.15in} & \hspace{-0.15in} & \hspace{-0.15in} &
\hspace{-0.15in} & \hspace{-0.15in} & \hspace{-0.15in} &
\hspace{-0.15in} &
\hspace{-0.15in} &  \\
\frac{-1}{2^{n-1}\times 3}\hspace{-0.15in} & \frac{1}{2^{n-1}\times 3}\hspace{%
-0.2in} & \frac{1}{2^{n-1}\times 3}\hspace{-0.15in} & \frac{1}{2^{n-1}\times 3%
}\hspace{-0.15in} & \frac{-1}{2^{n-1}\times 3}\hspace{-0.15in} & \frac{-1}{%
2^{n-1}\times 3}\hspace{-0.15in} & \frac{1}{2^{n-1}\times
3}\hspace{-0.15in} &
\frac{1}{2^{n}}\hspace{-0.15in} & \frac{1}{2^{n}}\hspace{-0.15in} & \frac{-1}{%
2^{n}}\hspace{-0.15in} & \frac{-1}{2^{n}}\hspace{-0.15in} &
0\hspace{-0.15in} &
.\hspace{-0.15in} & .\hspace{-0.15in} & \frac{1}{8}\hspace{-0.15in} & \frac{1}{8%
}\hspace{-0.15in} & \frac{-1}{8}\hspace{-0.15in} &
\frac{-1}{8}\hspace{-0.15in} & 0\hspace{-0.15in} &
\frac{1}{4}\hspace{-0.15in} & \frac{1}{4}\hspace{-0.15in}
& \frac{-1}{4}\hspace{-0.15in} & \frac{-1}{4}\hspace{-0.15in} & 0\hspace{-0.2in%
} & \frac{1}{2}\hspace{-0.15in} & \frac{1}{2}\hspace{-0.15in} & \frac{-1}{2}%
\hspace{-0.15in} & \frac{-1}{2}\hspace{-0.15in} & 0\hspace{-0.15in} & 0 \\
\hspace{-0.15in} & \hspace{-0.15in} & \hspace{-0.15in} &
\hspace{-0.15in} & \hspace{-0.15in} & \hspace{-0.15in} &
\hspace{-0.15in} & \hspace{-0.15in} & \hspace{-0.15in} &
\hspace{-0.15in} & \hspace{-0.15in} & \hspace{-0.15in} &
\hspace{-0.15in} & \hspace{-0.15in} & \hspace{-0.15in} &
\hspace{-0.15in} & \hspace{-0.15in} & \hspace{-0.15in} &
\hspace{-0.15in} & \hspace{-0.15in} & \hspace{-0.15in} &
\hspace{-0.15in} & \hspace{-0.15in} & \hspace{-0.15in} &
\hspace{-0.15in} & \hspace{-0.15in} & \hspace{-0.15in} &
\hspace{-0.15in} &
\hspace{-0.15in} &  \\
\frac{1}{3}\hspace{-0.15in} & \frac{-1}{3}\hspace{-0.15in} & \frac{2}{3}%
\hspace{-0.15in} & \frac{2}{3}\hspace{-0.15in} &
\frac{-2}{3}\hspace{-0.15in} & \frac{-2}{3}\hspace{-0.15in} &
\frac{-3}{4}\hspace{-0.15in} & 0\hspace{-0.15in} &
0\hspace{-0.15in} & 0\hspace{-0.15in} & 0\hspace{-0.15in} &
0\hspace{-0.15in} & .\hspace{-0.15in} & .\hspace{-0.15in} &
0\hspace{-0.15in} & 0\hspace{-0.15in} & 0\hspace{-0.15in} &
0\hspace{-0.15in} & 0\hspace{-0.15in} & 0\hspace{-0.15in} &
0\hspace{-0.15in} & 0\hspace{-0.15in} & 0\hspace{-0.15in} &
0\hspace{-0.15in} & 0\hspace{-0.15in} & 0\hspace{-0.15in} &
0\hspace{-0.15in} & 0\hspace{-0.15in}
& 0\hspace{-0.15in} & 1 \\
\hspace{-0.15in} & \hspace{-0.15in} & \hspace{-0.15in} &
\hspace{-0.15in} & \hspace{-0.15in} & \hspace{-0.15in} &
\hspace{-0.15in} & \hspace{-0.15in} & \hspace{-0.15in} &
\hspace{-0.15in} & \hspace{-0.15in} & \hspace{-0.15in} &
\hspace{-0.15in} & \hspace{-0.15in} & \hspace{-0.15in} &
\hspace{-0.15in} & \hspace{-0.15in} & \hspace{-0.15in} &
\hspace{-0.15in} & \hspace{-0.15in} & \hspace{-0.15in} &
\hspace{-0.15in} & \hspace{-0.15in} & \hspace{-0.15in} &
\hspace{-0.15in} & \hspace{-0.15in} & \hspace{-0.15in} &
\hspace{-0.15in} &
\hspace{-0.15in} &  \\
\frac{-(2^{n-1}-1)}{2^{n-1}\times 3}\hspace{-0.15in} & \frac{2^{n-1}-1}{%
2^{n-1}\times 3}\hspace{-0.15in} & \frac{-(2^{n}+1)}{2^{n-1}\times 3}\hspace{%
-0.2in} & \frac{-(2^{n}+1)}{2^{n-1}\times 3}\hspace{-0.15in} & \frac{2^{n}+1}{%
2^{n-1}\times 3}\hspace{-0.15in} & \frac{2^{n}+1}{2^{n-1}\times 3}\hspace{%
-0.2in} & \frac{2^{n+1}-1}{2^{n-1}\times 3}\hspace{-0.15in} & \frac{-1}{2^{n}}%
\hspace{-0.15in} & \frac{-1}{2^{n}}\hspace{-0.15in} & \frac{1}{2^{n}}\hspace{%
-0.2in} & \frac{1}{2^{n}}\hspace{-0.15in} & 0\hspace{-0.15in} & .\hspace{-0.2in%
} & .\hspace{-0.15in} & \frac{-1}{8}\hspace{-0.15in} & \frac{-1}{8}\hspace{%
-0.2in} & \frac{1}{8}\hspace{-0.15in} & \frac{1}{8}\hspace{-0.15in} & 0\hspace{%
-0.2in} & \frac{-1}{4}\hspace{-0.15in} & \frac{-1}{4}\hspace{-0.15in} & \frac{1%
}{4}\hspace{-0.15in} & \frac{1}{4}\hspace{-0.15in} & 0\hspace{-0.15in} & \frac{-%
1}{2}\hspace{-0.15in} & \frac{1}{2}\hspace{-0.15in} & \frac{1}{2}\hspace{-0.2in%
} & \frac{1}{2}\hspace{-0.15in} & 0\hspace{-0.15in} & -1%
\end{bmatrix}%
} $

This shows that $F_{n}$ is full. To show that it is the geodesic
between the points $\overset{\rightarrow }{a_{1}}$ and
$\overset{\rightarrow }{a}_{5n+3}$ in $F_{n},$ we show that any
proper subset \ $A$ of $F_{n}$ containing these two points is not
full. If possible suppose such a set $A$ is full. Then $A$ has to
contain the point $\overset{\rightarrow }{b}_{n}$ because no subset
of $D_{n},$ other than singleton, is full.

Let $k=|F_{n}|-|A|.$ As $A$ is full there exists atleast $k$ coordinates of
points of $F_{n}$ which donot occur as coordinates in the points of $A.$
(Because otherwise adding these $k$ points we get $F_{n}$ and we will be
adding less than $k$ coordinates. If $A$ is full then $F_{n}$ cannot be
good). Let $S$ denote these $k$ coordinates. The set $S$ cannot contain $%
x_{1},x_{2},x_{3},\alpha _{5n-1},\alpha _{5n-2},\alpha _{5n-8},y_{2}$ and $%
z_{3}$ as these are used by the points of $A.$ Among these $k$ coordinates
let $k_{i}$ be the number which are introduced in $A_{i},0\leq i\leq n.$ For
$i\geqslant 1,$ we have $0\leq k_{i}\leq 5$ and $k_{0}=0$ or $1.$ If $%
0<k_{i}<5$ for some $i\geqslant 1,$ (or if $k_{0}=1$ for $i=0$ ) the $k_{i}$
coordinates of $S$ introduced in $A_{i}$ are used in atleast $k_{i}+1$
points of $A_{i}$. So if $k_{0}=1$ or $0<k_{i}<5$ for some $i\geqslant 1,$
then more than $k$ points of $F_{n}$ cannot be in $A$ which is a
contradiction. In the case $k_{0}=0$ and $k_{i}=0$ or $k_{i}=5$ for $%
i\geqslant 1,$ clearly there exists an $i\geqslant 1$ with $k_{i}=5.$ But in
this case we have $k_{n-1}=k_{n}=0.$ If $k_{i}=5,$ then $A\cap A_{i}=\phi $
and if $k_{i}=0,$ then $A\cap A_{i}=A_{i}.$ Let $j$ be an index such that $%
k_{j}=5$ and $k_{j+1}=0.$ Then $A\cap A_{j+1}=A_{j+1}$ which is a
contradiction because $A_{j+1}$ uses coordinates introduced in $A_{j}$ which
are not used by points of $A.$ This shows $A$ is not full.

It can be seen that the $5$ rows of $M_{n}^{-1},$ from $(5m-1)th$ row to $\
(5m+3)$rd \ row, have row sums bounded by $C_{1}+C_{2}\sum_{i=1}^{m}\frac{1}{%
2^{i}}$ for some constants $C_{1}$ and $C_{2},$\ \ independent of $n.$ This
shows that as in higher dimensions, in the three dimensional case also
uniform boundedness of lengths of geodesics is not a necessary condition for
boundedness of solutions of $(1)$ for bounded function $f$. \bigskip

\textsl{Acknowledgement: } I thank Prof. M G Nadkarni for suggesting the
problem, fruitful discussions and encouragement.

\bigskip

\bigskip \textbf{References:}

\bigskip [1]\quad Cowsik R C, K\l opotowski A and Nadkarni M G, When is $%
f(x,y)=u(x)+v(y)$?, \textit{Proc. Indian Acad. Sci. (Math. Sci.)} 109 (1999)
57-64.

[2]\quad K\l opotowski A, Nadkarni M G and Bhaskara Rao K P S, When is $%
f(x_{1},x_{2},\ldots ,x_{n})=u_{1}(x_{1})+u_{2}(x_{2})+...+u_{n}(x_{n})$?,
\textit{Proc. Indian Acad. Sci. (Math. Sci.) } 113 (2003) 77-86.

[3]\quad K\l opotowski A, Nadkarni M G and Bhaskara Rao K P S, Geometry of
good sets in $n$-fold Cartesian products, \textit{Proc. Indian Acad. Sci.
(Math. Sci.)} 114 (2004) 181-197.

[4]\quad Nadkarni M G, Kolmogorov' s superposition theorem and sums of
algebras, \textit{The Journal of Analysis} vol. 12 (2004) 21-67.

[5]\quad Gowri Navada K, Some remarks on good sets, \textit{Proc. Indian
Acad. Sci. (Math. Sci.)} 114 (2004) 389-397.

[6] \ \ Gowri Navada K, Some further remarks on good sets, to appear in
\textit{Proc. Indian Acad. Sci. (Math. Sci.)}.

\end{document}